\documentclass[11pt, twoside]{article}
\usepackage{amssymb}
\usepackage{amsmath, amssymb}\usepackage{cite}\usepackage{mathrsfs}\usepackage[amsmath, thmmarks]{ntheorem}
\usepackage{listings}
\usepackage[titletoc]{appendix}\allowdisplaybreaks
\usepackage{multirow}

\makeatletter
\renewcommand\theequation{\thesection.\@arabic\c@equation}
\makeatother

\newtheorem{thm}{Theorem}[section]%
\newtheorem{lem}[thm]{Lemma}%

\newtheorem{Remark}[thm]{Remark}%
\newtheorem{Pro}[thm]{Proposition}%
\newtheorem{Con}[thm]{Conjecture}%
\newtheorem{Fac}[thm]{Fact}%
\newtheorem{Dec}[thm]{Declaration}%
\input cyracc.def

\def\f{\noindent}

\def\demo{\f{\bf Proof}\hskip10pt}

\def\qed{\hfill $\Box$}

\begin{document}

\title{{\bf A note on the horizontal class transposition group}}
\footnotetext{E-mail addresses: Junyao$_{-}$Pan@126.com}

\author{{\bf Junyao Pan}\\
{\footnotesize School of Cyber Science and Engineering, Wuxi University, Wuxi, Jiangsu, 214105, P. R. China}
}

\date{}
\maketitle

\noindent{\small {\bf Abstract:} Let $n$ be an integer with $n > 1$. For every $r$ satisfying the inequalities $0 \leq r < n$, the residue class modulo $n$ is defined as $r(n)=\{r + kn | k \in Z\}$, where $Z$ is the set of all integers. Then for $0 \leq r_1\neq r_2 < n$, the horizontal class transposition $\tau_{r_1(n), r_2(n)}$ is an involution that interchanges $r_1 + kn$ and $r_2 + kn$ for each integer $k$ and fixes everything else. The horizontal class transposition group $CT_n$ is generated by all horizontal class transposition $\tau_{r_1(n), r_2(n)}$. Let $N$ be the least common multiple of the numbers $2, 3, . . . , n$ and $CT_{(n)}=\langle CT_2,CT_3,...,CT_n\rangle$. In this note, we prove that for $n>3$, $CT_{(n)}\cong S_N$, where $S_N$ is the symmetric group of degree $N$. Thus, we solve a conjecture proposed by Bardakov and Iskra, which has been included in the kourovka notebook: Unsolved problems in group theory, Novosibirsk, 2026.

\vskip0.2cm
\noindent{\small {\bf Keywords:} Residue class; Horizontal class transposition; Symmetric group.


\section {Introduction}

Throughout this paper, $Z$ denotes the set of all integers; $S_n$ stands for the symmetric group of degree $n$; $n\mid m$ and $n\nmid m$ means that $n$ divides and does not divide $m$, respectively; $a\equiv b(mod~n)$ and $a\not\equiv b(mod~n)$ imply that $a$ is congruent and not congruent to $b$ modulo $n$, respectively.

Let $m$ be a natural number with $m > 1$. For every $r$ satisfying the inequalities $0 \leq r < m$, the \emph{residue class modulo $m$} is defined as $r(m) = r + mZ = \{r + km | k \in Z\}$.
If $r_1(m_1) \cap r_2(m_2) = \emptyset$, the \emph{class transposition} $\tau_{r_1(m_1), r_2(m_2)}$ is the product of independent transpositions, that is, $$\tau_{r_1(m_1), r_2(m_2)}=\prod_{k \in Z}(r_1+km_1,r_2+km_2).$$
If $m_1=m_2$ then $\tau_{r_1(m_1), r_2(m_2)}$ is called \emph{horizontal class transposition}; otherwise, $\tau_{r_1(m_1), r_2(m_2)}$ is called \emph{slanted class transposition}. S. Kohl \cite{K1} introduced the \emph{class transposition group} $CT(Z)$ that is generated by all class transpositions, and showed that this group has a number of remarkable properties, such as the $CT(Z)$ is simple, and contains all symmetric group $S_n$ and hence all finite groups. Additionally, Kohl proposed an open problem that is it true that there exists only a finite set of natural numbers that occur as the orders of the product of two class transpositions? In fact, some connections between this problem and the famous Collatz conjecture have been established in \cite{K2}, and further this open problem has been included in \cite[Problem 18.48]{KU}. Recently, Bardakov and Iskra made breakthrough progress on this problem, see references \cite{BA,BA1}.

Given an integer $n$ with $n > 1$. The \emph{horizontal class transposition group} $CT_n$ is generated by all horizontal class transpositions $\tau_{r_1(n), r_2(n)}$. Let $CT_{(n)}=\langle CT_2,CT_3,...,CT_n\rangle$ and $N$ be the least common multiple of the numbers $2, 3, . . . , n$. For the group $CT_{(n)}$, Bardakov and Iskra \cite{BA} proposed a conjecture that has been included in \cite[Problem 21.8]{KU}, as follows:

\begin{Con}\label{pan1-2}\normalfont([5, Problem 21.8])
Is it true that for $n > 3$ the group $CT_{(n)}$ is isomorphic to the symmetric group $S_N$.
\end{Con}
It is well-known that the symmetric groups are a class of important groups, and we are interested in the symmetric groups and their applications, refer to \cite{P1,P2,P3,P4}. So in this note, we investigate and resolve Conjecture\ \ref{pan1-2}. More precisely, we obtain the following theorem.

\begin{thm}\label{pan1-1}\normalfont
Let $CT_{(n)}=\langle CT_2,CT_3,...,CT_n\rangle$ and $N$ be the least common multiple of the numbers $2, 3, . . . , n$. Then for any $n > 3$, $CT_{(n)}\cong S_N$.
\end{thm}

\section {Preliminaries}

Firstly, we introduce some notions and notations that will be used. Let $n$ be a natural number with $n>1$. Then for any $k\in Z$, we set $R^k_n=\{kn,kn+1,...,kn+n-1\}$. Applying the definition of the horizontal class transposition, it follows that for any $g\in CT_n$, $(R^k_n)^g=R^k_n$ for each $k\in Z$, and further for $s\in R^0_n$, if $s^g=t$ then $t\in R^0_n$ and $(kn+s)^g=kn+t$ for each $k\in Z$. So we call that the $CT_n$ \emph{independently} and \emph{synchronously} acts on $R^k_n$ for all $k\in Z$. Let $CT^k_n$ be the group that is the restriction of $CT_n$ acting on $R^k_n$. Since $CT_n$ independently and synchronously acts on $R^k_n$ for all $k\in Z$, it follows that $CT_n\cong CT^k_n$ for any $k\in Z$. Moreover, let $R^{[k]}_n=R^0_n\cup R^1_n\cup\cdot\cdot\cdot\cup R^{k-1}_n$ and $CT^{[k]}_n$ denote the group that is the restriction of $CT_n$ acting on $R^{[k]}_n$, where $k\geq1$. Similarly, we have $CT_n\cong CT^{[k]}_n$. Since all transpositions generate the symmetric group, we deduce the following lemma immediately.

\begin{lem}\label{pan2-1}\normalfont
Let $n$ be a natural number with $n>1$. Then $CT_n\cong CT^k_n\cong S_n$ for any $k\in Z$, and $CT_n\cong CT^{[k]}_n\cong S_n$ for any $k>1$.
\end{lem}

\begin{lem}\label{pan2-2}\normalfont
Let $m$ and $n$ be two positive integers with $m=ln$ and $l>1$. Then $CT^{[l]}_n\leq CT^0_m$.
\end{lem}
\demo Note $R^0_m=R^{[l]}_n$. Then by Lemma\ \ref{pan2-1}, we deduce this lemma immediately.   \qed

\begin{lem}\label{pan2-22}\normalfont
Let $n$ be a natural number with $n>1$ and $(n+1)\nmid N$, where $N$ is the least common multiple of the numbers $2, 3, . . . , n$. Then $n+1=p^k$ for some prime $p$ and positive integer $k$. In addition, $Np$ is the least common multiple of the numbers $2, 3, . . . , n,n+1$.
\end{lem}
\demo If there exist two coprime positive integers $a,b$ such that $n+1=ab$ and $1<a,b<n+1$, then $a\mid N$, $b\mid N$ and $ab\mid N$, a contradiction. Thereby, $n+1=p^k$ for some $k\geq1$, where $p$ is a prime. Clearly, $p^{k-1}\mid N$ and so $Np$ is the least common multiple of the numbers $2, 3, . . . , n,n+1$. The proof of this lemma is completed.   \qed

Next we present some notions and notations about permutation groups, and for more details can be seen in \cite{D}. Let $G$ be a permutation group acting on a set $\Omega$, and $\Delta\subseteq\Omega$. Then the \emph{point stabilizer} of $\Delta$ in $G$ is denoted by $G_{(\Delta)}=\{g\in G|\alpha^g=\alpha,\forall\alpha\in\Delta\}$. In particular, if $\Delta=\{\alpha\}$ then $G_{(\Delta)}$ is abbreviated as $G_\alpha$. For a point $\alpha\in\Omega$, if $\Omega=\{\alpha^g|g\in G\}$ then $G$ is said to be \emph{transitive} on $\Omega$. If any $k$-tuple of distinct points can be mapped, by some element of $G$, to any other $k$-tuple of distinct points, then $G$ is said to be \emph{$k$-transitive} on $\Omega$. If $k>1$, then $G$ is $k$-transitive if and only if $G$ is $(k-1)$-transitive and $(k-1)$-points stabilizer is transitive on the set of remaining points. Next we state a well-known fact while is the key path to resolve Conjecture\ \ref{pan1-2}.

\begin{Fac}\label{pan1-3}\normalfont
Let $G$ be a $6$-transitive permutation group acting on a set $\Omega$. Then the group $G$ contains the alternating group on $\Omega$.
\end{Fac}

Now we provide the outline of our proof. Bardakov and Iskra \cite{BA} showed the Conjecture\ \ref{pan1-2} holds for $n=4$ by using GAP. We prove Conjecture\ \ref{pan1-2} by induction on $n$. So we assume that the Conjecture\ \ref{pan1-2} holds for $n$, where $n > 3$. If $(n+1)\mid N$, then $N$ is also the least common multiple of the numbers $2, 3, . . . , n,n+1$. Then by Lemma\ \ref{pan2-2} and inductive hypothesis, we infer that the Conjecture\ \ref{pan1-2} holds for $n+1$. So it suffices to consider the case $(n+1)\nmid N$. In this case, it follows from Lemma\ \ref{pan2-22} that $n+1=p^k$ for some prime $p$ and positive integer $k$. For convenience, we make the following statement and will not repeatedly emphasize it.

\begin{Dec}\label{pan2-222}\normalfont
In the following, $n$ always represents a positive inter such that $n>3$ and $n+1=p^k$ for some prime $p$ and positive integer $k$, $N$ always denotes the least common multiple of the numbers $2, 3, . . . , n$, and $M=\frac{N}{p^{k-1}}$.
\end{Dec}
One easily checks that $R^0_{Np}=R^{[p]}_N=R^{[M]}_{p^k}$. It follows from Lemma\ \ref{pan2-2} that $CT^{[p]}_N\leq CT^0_{Np}$ and $CT^{[M]}_{p^k}\leq CT^0_{Np}$. Then by Lemma\ \ref{pan2-1} we see that it suffices to prove that $\langle CT^{[p]}_N,CT^{[M]}_{p^k}\rangle=CT^0_{Np}$. Actually, we prove that $\langle CT^{[p]}_N,CT^{[M]}_{p^k}\rangle$ is a $6$-transitive permutation group on the set $R^0_{Np}$, and thus $\langle CT^{[p]}_N,CT^{[M]}_{p^k}\rangle$ contains the alternating group on the set $R^0_{Np}$. Moreover, the Lemma\ \ref{pan2-1} implies that $\langle CT^{[p]}_N,CT^{[M]}_{p^k}\rangle$ contains some odd permutations. Therefore, $\langle CT^{[p]}_N,CT^{[M]}_{p^k}\rangle=CT^0_{Np}\cong S_{Np}$. In order to prove that $\langle CT^{[p]}_N,CT^{[M]}_{p^k}\rangle$ is $6$-transitive, we introduce some notions and notations that will be used.

Let $R^i_N, R^{i+1}_N\subseteq R^{[p]}_N$. If $R^j_{p^k}\subseteq R^{[M]}_{p^k}$ such that $R^j_{p^k}\cap R^i_N\neq\emptyset$ and $R^j_{p^k}\cap R^{i+1}_N\neq\emptyset$, then $R^j_{p^k}$ is called a \emph{bridge} of $R^i_N$ and $R^{i+1}_N$. If $\big|R^j_{p^k}\cap R^i_N\big|=1~{\rm{and}}~\big|R^j_{p^k}\cap R^{i+1}_N\big|=p^k-1$, then $R^j_{p^k}$ is called a \emph{left single-log bridge} of $R^i_N$ and $R^{i+1}_N$; if $\big|R^j_{p^k}\cap R^i_N\big|=p^k-1~{\rm{and}}~\big|R^j_{p^k}\cap R^{i+1}_N\big|=1$, then $R^j_{p^k}$ is called a \emph{right single-log bridge} of $R^i_N$ and $R^{i+1}_N$. Indeed, $\langle CT^{[p]}_N,CT^{[M]}_{p^k}\rangle$ transitively acts on $R^0_{Np}$ by bridges, and so we next state some facts about bridges.

\begin{Pro}\label{pan2-3}\normalfont
For any $R^i_N, R^{i+1}_N\subseteq R^{[p]}_N$, there exists a bridge $R^j_{p^k}\subseteq R^{[M]}_{p^k}$ of $R^i_N$ and $R^{i+1}_N$.
\end{Pro}
\demo We prove this proposition by contradiction. Assume that there is no bridge of $R^i_N$ and $R^{i+1}_N$ for some $R^i_N, R^{i+1}_N\subseteq R^{[p]}_N$. Then we see that there exists a $R^j_{p^k}\subseteq R^{[M]}_{p^k}$ such that $R^{[i+1]}_N=R^{[j+1]}_{p^k}$. Thus, $Ni=jp^k\neq0$ and which imply $p\mid i$ and $M\mid j$. However, $0<i<p$ and $0<j<M$ indicate $p\nmid i$ and $M\nmid j$, a contradiction. The proof of this proposition is completed.   \qed

Applying the definitions of the left and right single-log bridge, we obtain the following remark.
\begin{Remark}\label{pan2-4}\normalfont
Let $R^j_{p^k}\subseteq R^{[M]}_{p^k}$ be a left single-log bridge of $R^i_N$ and $R^{i+1}_N$. Then $$R^j_{p^k}\cap R^i_N=\{jp^k\}=\{(i+1)N-1\}.$$ In this case, $jp^k$ is called a \emph{left single bridge pier}. Let $R^j_{p^k}\subseteq R^{[M]}_{p^k}$ be a right single-log bridge of $R^i_N$ and $R^{i+1}_N$. Then $$R^j_{p^k}\cap R^{i+1}_N=\{(j+1)p^k-1\}=\{(i+1)N\}.$$ In this case, $(i+1)N$ is called a \emph{right single bridge pier}.
\end{Remark}

\begin{Pro}\label{pan2-5}\normalfont
Let $\mathcal{B}$ be the set of all bridges of $R^i_N$ and $R^{i+1}_N$ for $i=0,1,...,p-2$. Then there exists at most one left single-log bridge and at most one right single-log bridge in $\mathcal{B}$.
\end{Pro}
\demo Assume that there exist two left single-log bridges in $\mathcal{B}$. Let $\big|R^j_{p^k}\cap R^i_N\big|=1$, $\big|R^j_{p^k}\cap R^{i+1}_N\big|=p^k-1$, $\big|R^{j'}_{p^k}\cap R^{i'}_N\big|=1$ and $\big|R^{j'}_{p^k}\cap R^{{i'}+1}_N\big|=p^k-1$. It follows from Remark\ \ref{pan2-4} that $(i+1)N-1=jp^k$ and $(i'+1)N-1=j'p^k$. Thus, $(i'-i)N=(j'-j)p^k$, however, $p^k\mid(i'-i)N$ and $N\mid(j'-j)p^k$ if and only if $i'=i$ and $j'=j$. Therefore, there exists at most one left single-log bridge in $\mathcal{B}$. Similarly, we can infer that there exists at most one right single-log bridge in $\mathcal{B}$. The proof of this proposition is completed.   \qed

Next we prove that $\langle CT^{[p]}_N,CT^{[M]}_{p^k}\rangle$ is transitive on the set $R^0_{Np}$ by means of bridges.

\begin{lem}\label{pan2-6}\normalfont
The group $\langle CT^{[p]}_N,CT^{[M]}_{p^k}\rangle$ acts transitively on the set $R^0_{Np}$.
\end{lem}
\demo Note that $R^0_{Np}=R^{[p]}_N=R^{[M]}_{p^k}=R^0_N\cup R^1_N\cup\cdot\cdot\cdot\cup R^{p-1}_N=R^0_{p^k}\cup R^1_{p^k}\cup\cdot\cdot\cdot\cup R^{M-1}_{p^k}$. Let $\alpha$ and $\beta$ be two distinct integers in $R^0_{Np}$. If $\alpha,\beta\in R^s_N$  for some $s\in R^0_p$, then by Lemma\ \ref{pan2-1} we see that there exists a $g\in CT^{[p]}_N$ such that $\alpha^g=\beta$. Similarly, if $\alpha,\beta\in R^t_{p^k}$ for some $t\in R^0_M$, then there exists a $g\in CT^{[M]}_{p^k}$ such that $\alpha^g=\beta$. Thereby, we assume that $\alpha$ and $\beta$ are in distinct sets of $\{R^0_N, R^1_N,...,R^{p-1}_N\}$. We claim that if $\alpha\in R^i_N$ and $\beta\in R^{i+1}_N$ for any $R^i_N, R^{i+1}_N\subseteq R^{[p]}_N$, then there exists a $g\in \langle CT^{[p]}_N,CT^{[M]}_{p^k}\rangle$ such that $\alpha^g=\beta$. It follows from Proposition\ \ref{pan2-3} that there exists a bridge $R^j_{p^k}\subseteq R^{[M]}_{p^k}$ such that $R^j_{p^k}\cap R^i_N\neq\emptyset$ and $R^j_{p^k}\cap R^{i+1}_N\neq\emptyset$. By Lemma\ \ref{pan2-1} we see that there exist $g_1\in CT^{[p]}_N$ and $g_2\in CT^{[M]}_{p^k}$ such that $\alpha^{g_1}\in R^j_{p^k}\cap R^i_N$ and $(\alpha^{g_1})^{g_2}\in R^{i+1}_N$. Applying Lemma\ \ref{pan2-1} again, we see that there exists a $g_3\in CT^{[p]}_N$ such that $[(\alpha^{g_1})^{g_2}]^{g_3}=\beta$. Let $g=g_1g_2g_3$. Clearly, $g\in \langle CT^{[p]}_N,CT^{[M]}_{p^k}\rangle$ such that $\alpha^g=\beta$, and so our claim is true. Let $\alpha\in R^s_N$ and $\beta\in R^{t}_N$ with $t-s>1$. Pick $\beta_i\in R^{s+i}_N$ where $i\in\{1,...,t-s\}$ and $\beta=\beta_{t-s}$. Then by applying our claim, we can deduce that there exists $g_1\in \langle CT^{[p]}_N,CT^{[M]}_{p^k}\rangle$ and $g_i\in \langle CT^{[p]}_N,CT^{[M]}_{p^k}\rangle$ such that $\alpha^{g_1}=\beta_1$ and $(\beta_{i-1})^{g_i}=\beta_i$ for $i=2,...,t-s$. Let $g=g_1g_2\cdot\cdot\cdot g_{t-s}$. It is obvious that $g\in \langle CT^{[p]}_N,CT^{[M]}_{p^k}\rangle$ such that $\alpha^g=\beta$. Moreover, it is clear that the status of $\alpha$ and $\beta$ is equivalent. Therefore, $\langle CT^{[p]}_N,CT^{[M]}_{p^k}\rangle$ is transitive on the set $R^0_{Np}$, as desired.  \qed

\section {Proof of Theorem 1.2}
Before proving the multiple transitivity, we provide some concepts and simple facts that will be used. Let $(R^{j}_{p^k})_2=(R^{j}_{p^k})\setminus\{jp^k+2\}$ for $j\in R^0_M$. For any two distinct $r_1,r_2\in(R^{0}_{p^k})_2$, we define the transposition $\omega_{r_1(p^k), r_2(p^k)}$ by the rule that $$\omega_{r_1(p^k), r_2(p^k)}=\prod^{M-1}_{i=0}(r_1+ip^k,r_2+ip^k).$$
By the definition of $CT^{[M]}_{p^k}$, we see that $(CT^{[M]}_{p^k})_2$ are generated by all $\omega_{r_1(p^k), r_2(p^k)}$, where $(CT^{[M]}_{p^k})_2$ denotes the stabilizer of $2$ in $CT^{[M]}_{p^k}$. Let $\mathcal{E}=\{2+ap^k{\big{|}}a\in M_E\}$, where $M_E=\{a\in R^0_M{\big{|}}~2\mid a\}$. Now we state some useful and simple facts about $(CT^{[M]}_{p^k})_2$ and $\mathcal{E}$.

\begin{itemize}
\item The $(CT^{[M]}_{p^k})_2$ is the symmetric group that independently and synchronously acts on the set $(R^{j}_{p^k})_2$ for $j\in R^0_M$;
\item For any $\Delta\subseteq\mathcal{E}$, the point stabilizer of $\Delta$ in $CT^{[M]}_{p^k}$ is still $(CT^{[M]}_{p^k})_2$;
\item Every point in $\mathcal{E}$ is neither a left single bridge pier nor a right single bridge pier;
\item $a\not\equiv b~(mod~N)$ for any $a,b\in\mathcal{E}$.
\end{itemize}
\begin{Pro}\label{pan2-77}\normalfont
There is no $a\in M_E$ such that $2+ap^k+bN$ is a left single bridge pier for some $b\in R^0_p$. Moreover, there exists at most one $b\in M_E$ such that $2+bp^k+cN$ is a right single bridge pier for some $c\in R^0_p$.
\end{Pro}
\demo Suppose that there is $a\in M_E$ such that $2+ap^k+bN$ is a left single bridge pier for some $b\in R^0_p$. Then by Remark\ \ref{pan2-4}, we have $2+ap^k+bN=xN-1$ for some $x\in R^0_p$. Thus, $3=(x-b)N-ap^k$, however, $2\mid a$ and $2\mid N$ indicate $2\mid 3$, a contradiction. It follows from Proposition\ \ref{pan2-5} that there exists at most one right single bridge pier. Suppose that there are two $b,b'\in M_E$ such that $2+bp^k+cN=2+b'p^k+dN$ is a right single bridge pier. Then we have $(b-b')p^k=(d-c)N$. It is easy to see that $b=b'$ and $c=d$, a contradiction. Thus, there exists at most one $b\in M_E$ such that $2+bp^k+cN$ is a right single bridge pier for some $c\in R^0_p$, as desired.  \qed

Let $\Delta\subseteq\mathcal{E}$ and $\Delta_N=\{x\in R^0_N\big|x\equiv a(mod~N),a\in\Delta\}$. Define $(R^{j}_{N})_{\Delta_N}=R^{j}_{N}\setminus(\Delta_N+jN)$ for $j\in R^0_p$, where $\Delta_N+jN=\{x+jN|x\in\Delta_N\}$. For any $s_1,s_2\in (R^{0}_{N})_{\Delta_N}$, we define the transposition $\omega_{s_1(N), s_2(N)}$ by the rule that $$\omega_{s_1(N), s_2(N)}=\prod^{p-1}_{j=0}(s_1+jN,s_2+jN).$$By the definition of $CT^{[p]}_{N}$, we see that $(CT^{[p]}_{N})_{\Delta_N}$ are generated by all $\omega_{s_1(N), s_2(N)}$, where $(CT^{[p]}_{N})_{\Delta_N}$ denotes the stabilizer of $\Delta_N$ in $CT^{[p]}_{N}$. Noting that the $(CT^{[p]}_{N})_{\Delta_N}$ is the symmetric group that independently and synchronously acts on the set $(R^{j}_{N})_{\Delta_N}$ for $j\in\ R^0_p$, and $(CT^{[p]}_{N})_{\Delta_N}=(CT^{[p]}_{N})_{\Delta}$. Next we give a key proposition.

\begin{Pro}\label{pan2-7}\normalfont
Let $\Delta\subseteq\mathcal{E}$. If there exist is no $2+bp^k\in \Delta$ such that $2+bp^k+cN$ is a right single bridge pier for some $c\in R^0_p$, then for any $(R^{i}_{N})_{\Delta_N}, (R^{i+1}_{N})_{\Delta_N}$, there exists a $(R^j_{p^k})_2$ such that $$(R^j_{p^k})_2\cap (R^{i}_{N})_{\Delta_N}\neq\emptyset~{\rm{and}}~ (R^j_{p^k})_2\cap (R^{i+1}_{N})_{\Delta_N}\neq\emptyset.$$
\end{Pro}
\demo It follows from Proposition\ \ref{pan2-3} that there is a bridge $R^j_{p^k}$ of $R^i_N$ and $R^{i+1}_N$. So we assume that $R^j_{p^k}\cap R^i_N=\{jp^k,1+jp^k,...,s+jp^k\}$ and $R^j_{p^k}\cap R^{i+1}_N=\{1+s+jp^k,2+s+jp^k,...,p^k-1+jp^k\}$ for some $s\in\{0,1,...,p^k-2\}$. Suppose that $R^j_{p^k}$ is neither a left single-log bridge nor a right single-log bridge. Then we have $1+jp^k\in(R^j_{p^k})_2\cap (R^{i}_{N})_{\Delta_N}$, and if $p=2$ then $p^k-1+jp^k\in(R^j_{p^k})_2\cap (R^{i+1}_{N})_{\Delta_N}$, and if $p$ is odd prime, then $p^k-2+jp^k\in(R^j_{p^k})_2\cap (R^{i+1}_{N})_{\Delta_N}$. Suppose that $R^j_{p^k}$ is a single-log bridge. Note that every point in $\mathcal{E}$ is neither a left single bridge pier nor a right single bridge pier. By Proposition\ \ref{pan2-77}, we see that there is no $2+ap^k\in \Delta$ such that $2+ap^k+bN$ is a left single bridge pier for some $b\in R^0_p$. Moreover, there exist is no $2+bp^k\in \Delta$ such that $2+bp^k+cN$ is a right single bridge pier for some $c\in R^0_p$. Thus, $jp^k\in(R^j_{p^k})_2\cap (R^{i}_{N})_{\Delta_N}$ {\rm{and}} $p^k-1+jp^k\in(R^j_{p^k})_2\cap (R^{i+1}_{N})_{\Delta_N}$, as desired.  \qed

\begin{lem}\label{pan2-8}\normalfont
Let $\Delta\subseteq\mathcal{E}$. If there exist is no $2+bp^k\in \Delta$ such that $2+bp^k+cN$ is a right single bridge pier for some $c\in R^0_p$, then the group $\langle (CT^{[p]}_{N})_{\Delta},(CT^{[M]}_{p^k})_2\rangle$ is transitive on the set $R^0_{Np}\setminus\Delta$.
\end{lem}
\demo Note that both $(CT^{[p]}_{N})_{\Delta}$ and $(CT^{[M]}_{p^k})_2$ are the symmetric groups that independently and synchronously acts on $(R^{j}_{N})_{\Delta_N}$ and $(R^{i}_{p^k})_2$ for $j\in\{0,1,...,p-1\}$ and $i\in\{0,1,...,M-1\}$, respectively. Moreover, the Proposition\ \ref{pan2-7} implies that for any $(R^{i}_{N})_{\Delta_N}, (R^{i+1}_{N})_{\Delta_N}$, there exists a $(R^j_{p^k})_2$ such that $$(R^j_{p^k})_2\cap (R^{i}_{N})_{\Delta_N}\neq\emptyset~{\rm{and}}~ (R^j_{p^k})_2\cap (R^{i+1}_{N})_{\Delta_N}\neq\emptyset.$$Proceeding as in the proof of Lemma\ \ref{pan2-6}, we can deduce this lemma immediately.  \qed

\begin{lem}\label{pan2-9}\normalfont
The group $\langle CT^{[p]}_N,CT^{[M]}_{p^k}\rangle$ is $6$-transitive on the set $R^0_{Np}$.
\end{lem}
\demo According to Lemma\ \ref{pan2-6}, we see that $\langle CT^{[p]}_N,CT^{[M]}_{p^k}\rangle$ is transitive on $R^0_{Np}$. Pick $\Delta=\{2\}$. Clearly, $2$ is not a right single bridge pier. It follows from Lemma\ \ref{pan2-8} that $\langle (CT^{[p]}_{N})_{\Delta},(CT^{[M]}_{p^k})_2\rangle$ is transitive on the set $R^0_{Np}\setminus\Delta$. Note that $\langle (CT^{[p]}_{N})_{\Delta},(CT^{[M]}_{p^k})_2\rangle$ is the subgroup of the stabilizer of $\Delta$ in $\langle CT^{[p]}_N,CT^{[M]}_{p^k}\rangle$. Therefore, $\langle CT^{[p]}_N,CT^{[M]}_{p^k}\rangle$ is $2$-transitive on the set $R^0_{Np}$.

Since $n>3$, we have $|\mathcal{E}|\geq6$. Then by Proposition\ \ref{pan2-77}, we see that there exists a subset $\Delta\in \mathcal{E}$ such that $|\Delta|=5$ and there exist is no $2+bp^k\in \Delta$ such that $2+bp^k+cN$ is a right single bridge pier for some $c\in R^0_p$. So we can continue the above process step by step, we deduce that $\langle CT^{[p]}_N,CT^{[M]}_{p^k}\rangle$ is $6$-transitive on the set $R^0_{Np}$.  \qed

Up to now we have completed the proof of Theorem\ \ref{pan1-1}.

\section*{Use of AI tools declaration}
The authors declare they have not used Artificial Intelligence (AI) tools in the creation of this article.

\section*{Conflict of interest}

The authors declare there is no conflicts of interest.




\begin{thebibliography}{99}

\bibitem{BA}
V. G. Bardakov and A. L. Iskra, Orders of products of horizontal class transpositions, \emph{Mat. Zametki}, 118:5 (2025), 654-669; \emph{Math. Notes}, 118:5 (2025), 921-932.

\bibitem{BA1}
V. G. Bardakov and A. L. Iskra, Orders of Products of Slanted Class Transpositions, arXiv:2504.08595.

\bibitem{D}
J.D. Dixon and B. Mortimer, {\it Permutation Groups}, Graduate Texts in Mathematics, {\bf vol. 163}, Springer, New York, 1996.

\bibitem{K1}
S. Kohl, A simple group generated by involutions interchanging residue classes of the integers, \emph{Math. Z.}, 264, no. 4 (2010), 927-938.

\bibitem{K2}
S. Kohl, The Collatz conjecture in a group theoretic context, \emph{J. Group Theory}, 20, no. 5 (2017), 1025-1030.

\bibitem{KU}
E. I. Khukhro and V. D. Mazurov, The kourovka notebook: Unsolved problems in group theory, Novosibirsk, 2026.


\bibitem{P1}
Junyao Pan,\ On a conjecture about orders of products of elements in the symmetric group, \emph{J. Pure. Appl. Algebra.}, 222, no. 2 (2018), 291-296.

\bibitem{P2}
Junyao Pan,\ On simultaneous conjugation of permutations, \emph{J. Pure. Appl. Algebra.}, 226, no. 1 (2022), Paper No.106804, 15 pp.

\bibitem{P3}
Junyao Pan,\ The full automorphism groups of general position graphs, {\it J. Combin. Theory Ser. A}, {\bf 201} (2024) 105800.

\bibitem{P4}
Junyao Pan,\ A note on the Cuntz algebra automorphisms, accepted to \emph{Israel J. Math.} (2025), DOI:10.1007/s11856-025-2851-x. Corpus ID: 275133798.


\end{thebibliography}
\end{document}